\documentclass[11pt]{amsart}

\usepackage{amsmath}
\usepackage{amssymb}
\usepackage[all]{xy}

\newtheorem{theorem}{Theorem}[section]
\newtheorem{claim}[theorem]{Claim}

\newtheorem{lemma}[theorem]{Lemma}

\theoremstyle{definition}
\newtheorem{definition}[theorem]{Definition}

\newtheorem{question}[theorem]{Question}

\theoremstyle{remark}
\newtheorem{remark}[theorem]{Remark}

\newcount\skewfactor
\def\mathunderaccent#1#2 {\let\theaccent#1\skewfactor#2
\mathpalette\putaccentunder}
\def\putaccentunder#1#2{\oalign{$#1#2$\crcr\hidewidth
\vbox to.2ex{\hbox{$#1\skew\skewfactor\theaccent{}$}\vss}\hidewidth}}

\def\smallbox#1{\leavevmode\thinspace\hbox{\vrule\vtop{\vbox
   {\hrule\kern1pt\hbox{\vphantom{\tt/}\thinspace{\tt#1}\thinspace}}
   \kern1pt\hrule}\vrule}\thinspace}

\newcommand{\cf}{{\rm cf}}

\def\qedref#1{$\qed_{\reforiginal{#1}}$}

\title[Magidor filters]{On the verge of inconsistency: \\
Magidor cardinals and Magidor filters}
\author{Shimon Garti}
\address{Institute of Mathematics,
 The Hebrew University of Jerusalem,
 Jerusalem 91904, Israel}
\email{shimon.garty@mail.huji.ac.il}

\author{Yair Hayut}
\address{Institute of Mathematics,
 The Hebrew University of Jerusalem,
 Jerusalem 91904, Israel}
\email{yair.hayut@mail.huji.ac.il}

\author{Saharon Shelah}
\address{Institute of Mathematics,
 The Hebrew University of Jerusalem,
 Jerusalem 91904, Israel,
 and Department of Mathematics
 Rutgers University
 New Brunswick, NJ 08854, USA}
\email{shelah@math.huji.ac.il}
\urladdr{http://www.math.rutgers.edu/\char`\~shelah}


\begin{document}
\let\labeloriginal\label
\let\reforiginal\ref

\begin{abstract}
We introduce a model-theoretic characterization of Magidor cardinals, from which we infer that Magidor filters are beyond ZFC-inconsistency.
\end{abstract}

\maketitle

\newpage

\section{introduction}

Magidor cardinals were defined in \cite{GaHa} through a combinatorial property of coloring functions. 
Recall that $\lambda\rightarrow[\lambda]^{\aleph_0\text{-bd}}_\lambda$ means that for every function $f$ from the bounded subsets of $\lambda$ of size $\aleph_0$ into $\lambda$ there exists a set $A\in[\lambda]^\lambda$ for which the restriction of $f$ to the bounded subsets of $A$ of size $\aleph_0$ omits at least one color.
The combinatorial definition reads as follows:

\begin{definition}
\label{mmm} Magidor cardinals. \newline
A cardinal $\lambda$ is a Magidor cardinal iff $\lambda\rightarrow[\lambda]^{\aleph_0\text{-bd}}_\lambda$.
\end{definition}

The term $\aleph_0\text{-bd}$ in the above definition means every countable set, regardless of its order type (though a similar concept of Magidority can be defined with an additional limitation on the order type).
The restriction to $\aleph_0$-bounded subsets in the domain of the coloring is obligatory if one wishes to accept the axiom of choice. 
If $\lambda\rightarrow[\lambda]^\omega_\lambda$ then $\lambda$ is called $\omega$-J\'onsson, and there are no such cardinals in ZFC as proved in  \cite{MR0209161}. 
Consequently, if $\lambda$ is a Magidor cardinal then $\lambda>\cf(\lambda)=\aleph_0$.

A primary question is what is the consistency strength of Magidor cardinals.
On the one hand, every Magidor cardinal is a J\'onsson cardinal. On the other hand, if $\lambda$ is I1 (or even I2) then $\lambda$ is a Magidor cardinal. There is some evidence that the consistency strength of Magidor cardinals is far from the parallel strength of J\'onsson cardinals, see \cite{GaHa}. In the current paper we give another evidence, through the notion of Magidor filters. Recall that a filter $F$ over $\lambda$ is uniform iff every element of $F$ is of size $\lambda$.

\begin{definition}
\label{mfdef} Magidor filters. \newline 
Let $F$ be a filter over $\lambda$. \newline 
We say that $F$ is a Magidor filter iff
\begin{enumerate}
\item [$(\aleph)$] $F$ is a uniform filter.
\item [$(\beth)$] $F$ contains all the end-segments of $\lambda$.
\item [$(\gimel)$] For every coloring $c:[\lambda]^{\aleph_0\text{-bd}}\rightarrow\lambda$ there exists an element $x\in F$ for which $c''[x]^{\aleph_0\text{-bd}}\neq\lambda$.
\end{enumerate}
\end{definition}

The uniformity requirement is crucial here, as otherwise the concept of Magidor filter would be meaningless. 
The requirement that each end-segment belongs to $F$ can be satisfied by extending $F$, as every uniform extension of a Magidor filter preserves the Magidority.

The main theorem of this paper is that there are no Magidor filters, assuming the axiom of choice. Recall that every measurable cardinal carries a J\'onsson filter. Moreover, every singular cardinal limit of measurable cardinals carries a J\'onsson filter by the proof of Prikry, \cite{MR0262075}. So the main result here for Magidority stands in sharp contrast with J\'onssonicity, and demonstrates in another way the discrepancy between the notions of Magidority and J\'onssonicity.

Our notation is standard, and elaborated in \cite{GaHa}. For a general background we suggest \cite{MR1994835}. We sketch several conventions to be used in the present paper.
By the existence of an elementary embedding we mean that there is some transitive model $M$ of ZFC and a non-trivial elementary embedding $\jmath$ from $M$ into some other transitive model of ZFC. 
The model $M$ need not be a proper class, and occasionally it will not be a model of all ZFC axioms. In which cases we assume that $M$ satisfies enough ZFC in order to carry out the relevant proof.
Assuming the axiom of choice, the non-triviality yields an ordinal $\beta$ so that $\jmath(\beta)>\beta$. The first such ordinal is called the critical point of $\jmath$ and denoted ${\rm crit}(\jmath)$.
An ordinal $\eta$ is a fixed point of $\jmath$ iff $\jmath(\eta)=\eta$. A fixed point $\eta$ is non-trivial if $\eta>{\rm crit}(\jmath)$.

If $g:\nu\rightarrow\mu$ and $x\subseteq\mu$ then $g^{-1}[x]=\{\alpha\in\nu: g(\alpha)\in x\}$.
The arrows notation with $\aleph_0\text{-bd}$ as a superscript is coherent with the common arrows notation. If $F$ is a filter over $\lambda$ then $\lambda\rightarrow[F]^{\aleph_0\text{-bd}}_\theta$ means that the set from which we omit colors is an element of $F$. If $\lambda\rightarrow[F]^{\aleph_0\text{-bd}}_\lambda$ then there is an ordinal $\alpha<\lambda$ such that $\lambda\rightarrow[F]^{\aleph_0\text{-bd}}_{\alpha, <\alpha}$. 
This notation means that the number of attained colors is less than $\alpha$.
We denote the first such ordinal by $\alpha_M(F)$. 

A word about the axiom of choice is in order. The proof of the non-existence of Magidor filters is carried out in ZFC. If we drop AC (even upon replacing it by some weak versions of choice) then Magidor filters may exist. 
Moreover, they form a natural niche of large cardinals in ZF.
In some sense, Magidor filters continue the table of large cardinals beyond $\omega$-J\'onssonicity, since the existence of such a filter implies the existence of an unbounded set of $\omega$-J\'onsson cardinals. The fact that the concept of Magidor filters comes from the ZFC acceptable notion of Magidor cardinals in a fairly natural way, enables us to go further afield and gives another aspect to the interplay between large cardinals and the axiom of choice.

We thank the referee for the neat and rapid work on this paper.

\newpage

\section{A model-theoretic characterization of Magidor cardinals}

Most of large cardinals, including J\'onsson cardinals, have several natural definitions. This fact proves fruitful in many respects. For the main result we need a model theoretic characterization of Magidor cardinals. For a similar characterization of J\'onsson cardinals see Tryba, \cite{MR788256}.
Ahead of proving our characterization we need to show that each Magidor cardinal $\lambda$ is $\omega$-closed, i.e. $\alpha<\lambda\Rightarrow\alpha^\omega<\lambda$. This is the content of the following:

\begin{claim}
\label{oomegaclosed} Every Magidor cardinal is $\omega$-closed.
\end{claim}

\par\noindent\emph{Proof}. \newline 
Let $\lambda$ be a Magidor cardinal, and assume toward contradiction that $\mu<\lambda$ while $\mu^\omega\geq\lambda$. Since $\cf(\lambda)=\omega\neq \cf(\mu^\omega)$ we know that $\mu^\omega>\lambda$, though this fact has no significant role in the proof. From the limitude of $\lambda$ we have $\mu^{+3}<\lambda$.

We choose a sequence of sets $\langle C_i:i<\lambda\rangle$ with the following three properties:
\begin{enumerate}
\item [$(\alpha)$] Each $C_i$ is a bounded subset of $\lambda$, closed in its supremum.
\item [$(\beta)$] The order type of $C_i$ is $\mu$, for every $i<\lambda$.
\item [$(\gamma)$] For every $\delta<\lambda$ so that $\cf(\delta)=\mu^{+3}$ and for every club $E\subseteq\delta$ we can find an ordinal $i<\lambda$ for which $C_i\subseteq E$.
\end{enumerate}
The existence of such a sequence follows from the club guessing theorems of pcf, see \cite{MR1318912} or \cite{MR2768693}. Let $D$ be $\{u\subseteq\lambda: {\rm otp}(u)=\omega, {\rm sup}(u)<\lambda\}$, so $D$ is just the collection of bounded subsets of $\lambda$ of order type $\omega$. 
Let $S^\lambda_\omega$ be the set $\{\alpha<\lambda:\cf(\alpha)=\omega\}$.
We choose a function $d:D\rightarrow \lambda$ such that ${\rm Rang}(d\upharpoonright \{u\in D:u\subseteq C_i\cap S^\lambda_\omega\})=\lambda$ for every $i<\lambda$. The existence of $d$ follows from the fact that we have but $\lambda$-many $C_i$-s and $|\{u\subseteq C_i\cap S^\lambda_\omega:{\rm otp}(u)=\omega\}|>\lambda$ for every $i<\lambda$, hence we can choose the values of $d$ by induction on $i<\lambda$.

We shall define a coloring $c$ which exemplifies the non-Magidority of $\lambda$. Without loss of generality, ${\rm dom}(c)=\{u\subseteq\lambda: {\rm sup}(u)<\lambda, {\rm otp}(u)=\omega\cdot\omega\}$, since we can define $c(u)=0$ whenever $u\in[\lambda]^{\aleph_0\text{-bd}}$ and ${\rm otp}(u)\neq\omega\cdot\omega$.
If $u\in{\rm dom}(c)$, let $\langle\alpha_{u,i}:i<\omega\cdot\omega\rangle$ be an increasing enumeration of the members of $u$, and set:
$$
c(u)=d(\{\alpha:\exists n, \alpha=\bigcup\limits_{\ell\in\omega}\alpha_{u,\omega\cdot n+\ell}\}).
$$
Assume now that $A\subseteq\lambda$, the order type of $A$ is $\mu^{+3}$ and $A$ is bounded in $\lambda$. Let $B(A)=\{c(u):u\subseteq A, u\in{\rm dom}(c)\}$. We shall prove that $B(A)=\lambda$ for every $A$ as above, thus arriving at a contradiction. Indeed, given $S\in[\lambda]^\lambda$ let $A$ be the first $\mu^{+3}$ members of $S$. Since $B(A)=\lambda$ we conclude that $c$ omits no colors on $[S]^{\aleph_0\text{-bd}}$, contradicting the Magidority of $\lambda$.

Given any set $A$, we denote the closure of $A$ under the order topology by $c\ell(A)$.
Fix $A$ and let $\langle\alpha_i:i\leq\mu^{+3}\rangle$ enumerate the members of $c\ell(A)$, stipulating $\delta=\alpha_{\mu^{+3}}$. Choose any $\gamma<\lambda$. We shall prove that $\gamma\in B(A)$. Notice that $\cf(\delta)=\mu^{+3}$, and $\alpha_{i+1}\in A$ for every $i<\mu^{+3}$. Set $E=\{\alpha_{\omega\cdot i}: i<\mu^{+3}\}$, so $E$ is a club subset of $\delta$. By the choice of the club guessing sequence, there exists an ordinal $i_*<\lambda$ for which $C_{i_*}\subseteq E$.

By the definition of $d$ we can choose an increasing $\omega$-sequence $\langle\beta_n:n\in\omega\rangle$ of ordinals from $C_{i_*}\cap S^\lambda_\omega$ so that $d(\{\beta_n:n\in\omega\})=\gamma$. Recall that $C_{i_*}\subseteq E$, so each $\beta_n$ has the form $\alpha_{\omega\cdot i(n)}$ for some $i(n)<\mu^{+3}$, and $i(n)<i(n+1)$ for every $n\in\omega$.
Likewise, $\cf(\beta_n)=\omega$ for every $n\in\omega$.

We wish to define a set of order type $\omega\cdot\omega$ out of the $\beta_n$-s, with the goal being to define a set $u$ for which $c(u)=\gamma$.
For every $n\in\omega$ we choose an increasing sequence of ordinals $\langle j_{n,\ell}:\ell\in\omega\rangle$ such that:
\begin{enumerate}
\item [$(\alpha)$] $\bigcup\limits_{\ell\in\omega}j_{n,\ell}= \alpha_{\omega\cdot i(n)}$ for every $n\in\omega$.
\item [$(\beta)$] $j_{n,\ell}$ is a successor ordinal for every $n,\ell\in\omega$.
\item [$(\gamma)$] ${\rm max}\{\alpha_{\omega\cdot i(m)}:m<n\}< j_{n,0}$ for every $n\in\omega$.
\end{enumerate}

Let $u=\{\alpha_{j_{n,\ell}}:n,\ell\in\omega\}$. By $(\beta)$ we have $u\subseteq A$, so $u\subseteq\delta$ as $A\subseteq \delta$. By the construction, ${\rm otp}(u)=\omega\cdot\omega$. Observe that $\alpha_{\omega\cdot i(n)}={\rm sup}\{\beta\in u:{\rm otp}(u\cap\beta)< \omega\cdot n\}$ for every $n\in\omega$ (again, by the construction), and hence $c(u)=d(\{\alpha_{\omega\cdot i(n)}:n\in\omega\})=d(\{\beta_n: n\in\omega\})=\gamma$, so the proof is accomplished.

\hfill \qedref{oomegaclosed}

\begin{theorem}
\label{t1} Magidority and elementary embeddings. \newline 
Assume $\lambda>\cf(\lambda)=\omega$. \newline
The following are equivalent:
\begin{enumerate}
\item [$(a)$] $\lambda$ is a Magidor cardinal.
\item [$(b)$] For every $\gamma>\lambda$ there is a triple $(M,\bar{M},\jmath)$ such that $M\subseteq{\rm V}$ is a transitive set for which $\lambda+1\subseteq M$, $\jmath:M\rightarrow{\rm V}_\gamma$ an elementary embedding such that ${\rm crit}(\jmath)<\lambda=\jmath(\lambda)$, $\bar{M}=(M_n:n\in\omega)$, $M=\bigcup\limits_{n\in\omega}M_n$ and $M_n\prec M$ for every $n\in\omega$, for every $\delta<\lambda$ there exists $n_\delta\in\omega$ such that $n\in[n_\delta,\omega)\Rightarrow [M_n\cap\delta]^{\aleph_0}\subseteq M$ and finally $\lambda\nsubseteq\jmath'' M$.
\item [$(c)$] There exists some ordinal $\gamma>\lambda$ for which there is a triple $(M,\bar{M},\jmath)$ such that $M\subseteq{\rm V}$ is a transitive set for which $\lambda+1\subseteq M$, $\jmath:M\rightarrow{\rm V}_\gamma$ an elementary embedding such that ${\rm crit}(\jmath)<\lambda=\jmath(\lambda)$, $\bar{M}=(M_n:n\in\omega)$, $M=\bigcup\limits_{n\in\omega}M_n$ and $M_n\prec M$ for every $n\in\omega$, for every $\delta<\lambda$ there exists $n_\delta\in\omega$ such that $n\in[n_\delta,\omega)\Rightarrow [M_n\cap\delta]^{\aleph_0}\subseteq M$ and finally $\lambda\nsubseteq\jmath'' M$.
\end{enumerate}
\end{theorem}

\par\noindent\emph{Proof}. \newline 
For the direction $(a)\Rightarrow (b)$ fix any ordinal $\gamma>\lambda$. Let $\mathfrak{B}$ be an expansion of $({\rm V}_\gamma,\in)$ by Skolem functions. For every set $S\subseteq\lambda$ let $c\ell_{\mathfrak{B}}(S)$ be a minimal elementary submodel $M$ of $\mathfrak{B}$ for which $S\subseteq M$ and $[M\cap\lambda]^{\aleph_0}\subseteq M$. Such a model $M$ can be created by induction on $\omega_1$, using the L\"owenheim-Skolem theorem and adding the $\omega$-subsets at each step. Notice that $|c\ell_{\mathfrak{B}}(S)|\leq |S|^{\aleph_0}\cdot 2^{\aleph_0}$. In particular, if $u\in[\lambda]^{\aleph_0}$ then $|c\ell_{\mathfrak{B}}(u)|=2^{\aleph_0}$.

For every $u\in[\lambda]^{\aleph_0}$ fix an enumeration $\langle a_{u,\zeta}: \zeta<2^{\aleph_0}\rangle$ of the members of $c\ell_{\mathfrak{B}}(u)$. Let $\mathcal{A}\subseteq[\lambda]^{\aleph_0\text{-bd}}$ be a maximal almost disjoint family such that $a\in\mathcal{A}\Rightarrow{\rm otp}(a)=\omega$. Denote $\bigcup\{c\ell_{\mathfrak{B}}(u):u\in[\lambda]^{\aleph_0\text{-bd}}\}$ by $T$.

We define two functions, $c_0$ and $c_1$. The function $c_0:[\lambda]^{\aleph_0\text{-bd}}\rightarrow 2^{\aleph_0}\times\omega_1$ would give a pair of ordinals for the coloring $c_1$, and $c_1$ would be the coloring for which we employ assumption $(a)$. Our demand from $c_0$ is that for every $a\in\mathcal{A}$ the following statement will be satisfied:
$$
\forall b\in[a]^{\aleph_0}, {\rm Rang}(c_0\upharpoonright [b]^{\aleph_0})=2^{\aleph_0}\times\omega_1.
$$
Now we define $c_1:[\lambda]^{\aleph_0\text{-bd}}\rightarrow T$ as follows. If $u\in[\lambda]^{\aleph_0\text{-bd}}, {\rm otp}(u)=\delta+\omega$ for some limit ordinal $\delta$, $\{\alpha_{u,i}:i<\delta+\omega\}$ an increasing enumeration of the members of $u$, $v_u=\{\alpha_{u,\delta+n}:n\in\omega\}$ is the upper part of $u$, $c_0(v_u)=(\zeta,\varepsilon)$ and $a_{\{\alpha_{u,i}:i<\delta\cap\varepsilon\},\zeta}$ is an ordinal then we let $c_1(u)=a_{\{\alpha_{u,i}:i<\delta\cap\varepsilon\},\zeta}$. In all other cases, $c_1(u)=0$.
Since $|T|=\lambda$ we can use the Magidority assumption $(a)$ and choose a set $A\in[\lambda]^\lambda$ and an ordinal $\tau\in T$ such that $\forall u\in[A]^{\aleph_0\text{-bd}}, c_1(u)\neq\tau$.

Let $\langle\mu_n:n\in\omega\rangle$ be an increasing sequence of cardinals such that $\lambda=\bigcup\limits_{n\in\omega}\mu_n$. For every $n\in\omega$ let $M_n^*=c\ell_{\mathfrak{B}}(A\cap\mu_n)$. As noted above, $|M_n^*|\leq\mu_n^{\aleph_0}<\lambda$. By virtue of Skolem (i.e., by adding the Skolem functions to $\mathfrak{B}$ and taking an elementary submodel of it) $M_n^*\prec\mathfrak{B}$ and since $M_n^*\subseteq M_{n+1}^*$ we infer that $M_n^*\prec M_{n+1}^*$ for every $n\in\omega$. Let $M^*$ be $\bigcup\limits_{n\in\omega}M_n^*$ and let $M$ be the Mostowski collapse of $M^*$ by $\pi$.

The key-point in this part of the proof is that $\tau\notin M^*$. For showing this fact, assume towards contradiction that $\tau\in M^*$, and choose some $n\in\omega$ for which $\tau\in M_n^*$, i.e. $\tau\in c\ell_{\mathfrak{B}}(A\cap\mu_n)$. Pick up a set $u\in[A\cap\mu_n]^{\aleph_0}$ such that $\tau\in c\ell_{\mathfrak{B}}(u)$. By the fixed enumerations of $c\ell_{\mathfrak{B}}(u)$ we may write $\tau=a_{u,\zeta}$ for some $\zeta<2^{\aleph_0}$. Denote the order-type of $u$ by $\varepsilon$.

Choose some $a_0\subseteq A\setminus\mu_n, {\rm otp}(a_0)=\omega$ and $a_0$ is bounded in $\lambda$. Let $u_0=u\cup a_0$ so ${\rm otp}(u_0)=\varepsilon+\omega$. Let $n_0$ be the first natural number so that $u_0\subseteq\mu_{n_0}$. We choose now $a_1\in[A\setminus\mu_{n_0}]^{\aleph_0\text{-bd}}$ and $a_2\in\mathcal{A}$ such that $b=a_1\cap a_2$ is infinite. By the above properties of $c_0$, there exists $v_0\in[b]^{\aleph_0}$ such that $c_0(v_0)=(\zeta,\varepsilon)$.

Define $u_1=u_0\cup v_0$ and observe that ${\rm otp}(u_1)=\varepsilon+\omega+\omega$. Let $\delta=\varepsilon+\omega$, so $\delta$ is a limit ordinal. We claim that $c_1(u_1)=\tau$ which contradicts the choice of $\tau$. Indeed, recall that ${\rm otp}(u_1)=\delta+\omega$ and $\delta$ is the limit ordinal $\varepsilon+\omega$. By definition, $c_1(u_1)=a_{\{\alpha_{u_1,i}:i<\delta\cap\varepsilon\},\zeta}$, as $c_0(\{\alpha_{u,\delta_n}:n\in\omega\})=c_0(v_0)=(\zeta,\varepsilon)$. However, $\{\alpha_{u_1,i}:i<\delta\cap\varepsilon\}=\{\alpha_{u_1,i}:i<\varepsilon\}=u$ and hence $c_1(u_1)=a_{u,\zeta}=\tau$, a contradiction.

Let $\jmath$ be the inverse $\pi^{-1}$ of the Mostowski collapse. Notice that $\jmath$ is an elementary embedding from $M$ into ${\rm V}_\gamma$.
Observe also that $\jmath(\tau)\neq\tau$ since $\tau\notin\jmath'' M$, and hence $\jmath$ has a critical point below $\lambda$ and $\jmath(\lambda)=\lambda$. Let $M_n$ be $\pi(M_n^*)$ for every $n\in\omega$, and verify that all the requirements in $(b)$ are satisfied.

As $(b)\Rightarrow(c)$ trivially, we are left with $(c)\Rightarrow(a)$. Choose a quadruple $(\gamma,M,\bar{M},\jmath)$ as alleged in $(c)$, and assume toward contradiction that $\lambda\nrightarrow[\lambda]^{\aleph_0\text{-bd}}_\lambda$. Let $c:[\lambda]^{\aleph_0\text{-bd}}\rightarrow\lambda$ exemplify this fact. Notice that $c\in {\rm V}_\gamma$ and since $\jmath'' M\prec{\rm V}_\gamma$ we may assume without loss of generality that $c\in \jmath'' M$.

Fix an ordinal $\eta\in\lambda\setminus\jmath'' M$ (recall that $\lambda\nsubseteq\jmath'' M$) and an increasing sequence of regular cardinals $\langle\mu_n:n\in\omega\rangle$ which converges to $\lambda$ and begins with $\mu_0=\aleph_0$. By induction on $n\in\omega$ we choose a pair $(\ell_n,A_n)$ such that:
\begin{enumerate}
\item [$(\aleph)$] $m<n<\omega\Rightarrow \ell_m<\ell_n$.
\item [$(\beth)$] $A_n\subseteq \jmath'' M_{\ell_n}\cap\mu_{\ell_n}$ and $|A_n|=\mu_n$ for every $n\in\omega$.
\item [$(\gimel)$] $m<n<\omega\Rightarrow A_n\cap\mu_{\ell_m}=\emptyset$.
\end{enumerate}
How do we choose them? For $n=0$ let $\ell_0=0$ and $A_0=\{p\in\omega:p$ is a prime number$\}$. At the stage of $n+1$ we choose $\ell>\ell_n$ such that $|\jmath'' M\cap\mu_\ell|>\mu_{n+1}$ (such an $\ell$ exists as $|\jmath'' M\cap\lambda|=\lambda$). Recall that $M=\bigcup\limits_{n\in\omega}M_n$ so $\jmath'' M=\bigcup\limits_{n\in\omega}\jmath'' M_n$, and choose $\ell_{n+1}>\ell$ so that $|\jmath'' M_{\ell_{n+1}}\cap\mu_{\ell_{n+1}}|>\mu_{n+1}$. Consequently, there is a set $A_{n+1}\subseteq \jmath'' M_{\ell_{n+1}}\cap\mu_{\ell_{n+1}}$ so that $|A_{n+1}|=\mu_{n+1}$. By chopping an initial segment of $A_{n+1}$ we may assume that $A_{n+1}\cap\mu_{\ell_n}=\emptyset$, thus accomplishing the construction.

Set $A=\bigcup\limits_{n\in\omega}A_n$, so $A\in[\lambda]^\lambda$. Let us show that $\eta\notin c\upharpoonright[A]^{\aleph_0\text{-bd}}$. Indeed, let $u\subseteq A$ be any bounded set of size $\aleph_0$. By $(\gimel)$ above there is $n_0\in\omega$ such that $u\subseteq\bigcup\limits_{n\leq n_0}A_n$. By assumption $(c)$ of our theorem there is $n_1\in\omega$ such that $u\in\jmath'' M_{n_1}$, and hence $u\in\jmath'' M$. However, $c\in\jmath'' M$ as well and hence $c(u)\in\jmath'' M$ so $\eta\neq c(u)$ by the choice of $\eta$, a contradiction to the defining property of the coloring $c$.

\hfill \qedref{t1}

\begin{remark}
\label{rrremark} In the above theorem, one can show that for every $\delta<\lambda$ and every large enough $n\in\omega$ we have $[M_n\cap\delta]^{\aleph_0}\subseteq M_n$. Moreover, a similar statement for unbounded $\omega$-sequences of $M_n$ follows from the above proof. The omitted colors in the proof come from bounded $\omega$-sequences which are contained in $M$ but in none of the $M_n$-s.
\end{remark}

\hfill \qedref{rrremark}

The proof is suggestive in another direction. Usually, given a large enough cardinal we are confronted with the opposite situation in which the domain of the embedding is ${\rm V}$ (or some portion ${\rm V}_\gamma$ of it) and the range is some transitive model $M$. The more we require from $M$, the larger our cardinal is. In particular, one can step up in the chart of large cardinals by asking for strong closure properties from the $M$ side. The same holds in our characterization, but we have a glass ceiling. Recall that in ZFC there are no $\omega$-J\'onsson cardinals, by \cite{MR0311478}.

\begin{claim}
\label{cor1} $\omega$-J\'onsson cardinals. \newline 
Assume $\lambda$ is an infinite cardinal. \newline 
There is no triple $(\gamma,M,\jmath)$ such that:
\begin{enumerate}
\item [$(a)$] $\gamma>\lambda+1$.
\item [$(b)$] $M\subseteq{\rm V}$ is a transitive set.
\item [$(c)$] $\jmath:M\rightarrow{\rm V}_\gamma$ is an elementray embedding.
\item [$(d)$] $\lambda+1\subseteq M$ but $\lambda\nsubseteq\jmath'' M$.
\item [$(e)$] ${\rm crit}(\jmath)<\lambda=\jmath(\lambda)$.
\item [$(f)$] $[M\cap\lambda]^{\aleph_0}\subseteq M$.
\end{enumerate}
\end{claim}

\par\noindent\emph{Proof}. \newline 
By \cite{MR0209161}, $\lambda$ is not $\omega$-J\'onsson in ${\rm V}$, so also not in $M$ by elementarity. Choose a function $c:[\lambda]^\omega\rightarrow\lambda$ which exemplifies this fact. Without loss of generality, $c\in\jmath'' M$. We indicate that one may choose (by elementarity) $d\in M$ which exemplifies this property of $\lambda$ in $M$, and then define $c=\jmath'' d\in\jmath'' M$.

By our assumptions $\lambda\nsubseteq\jmath'' M$ so we can choose an ordinal $\eta\in\lambda\setminus\jmath'' M$. Let $A$ be the set $\jmath'' M\cap\lambda$, so $A\in[\lambda]^\lambda$. Let $u\in[A]^{\aleph_0}$ be any set. Since $u\subseteq A=\jmath'' M\cap\lambda$ we infer that $u\in\jmath'' M$.
Indeed, enumerate the members of $u$ by $\{u_n:n\in\omega\}$. Since $u\subseteq \jmath'' M$ we can choose $v_n\in M$ for every $n\in\omega$ so that $\forall n\in\omega,\jmath(v_n)=u_n$. Let $v=\{v_n:n\in\omega\}$, and notice that $\jmath(v)=u$. However, $v\in M$ by assumption $(f)$, so $u\in\jmath'' M$ and hence $c(u)\in\jmath'' M$. Consequently, $\eta\neq c(u)$ for every such $u$, a contradiction.

\hfill \qedref{cor1}

\newpage 

\section{Magidor filters}

For proving the main theorem, we need some additional facts. Recall that a filter $F$ over $\lambda$ is not $(\omega,\theta)$-regular iff for every $A\subseteq F,|A|=\theta$ one can find $B\subseteq A,|B|=\aleph_0$ such that $\bigcap\{x:x\in B\}$ is not empty. J\'onsson filters (and hence Magidor filters) possess some degree of irregularity. The proof of the following lemma is derived from the ideas of Prikry, \cite{MR0262075}.

\begin{lemma}
\label{lem1} Let $F$ be a Magidor filter over $\lambda$. \newline 
Then $F$ is not $(\omega,\theta)$-regular for some $\theta=\cf(\theta)<\lambda$. Actually, $\alpha_M(F)$ can serve as $\theta$.
\end{lemma}

\par\noindent\emph{Proof}. \newline 
If $F$ is a Magidor filter then, in particular, $F$ is a J\'onsson filter. One has to convert colorings of $\omega$-bounded subsets into coloring of finite subsets (this is done, e.g., in \cite{GaHa} with respect to being a Magidor cardinal). By Lemma 1.40 in \cite{MR0262075} it follows that $F$ is not $(\omega,\theta)$-regular for some $\theta=\cf(\theta)<\lambda$ (see also Tryba, \cite{MR877853}). 

\hfill \qedref{lem1}

\begin{remark}
\label{r1} Actually, one can prove stronger irregularity properties for J\'onsson and Magidor filters. The proof can be extracted from \cite{MR0262075}. Irregularity gives some kind of reflection, to be used in the main theorem.
\end{remark}

\hfill \qedref{r1}

If $\lambda$ carries a Magidor filter $F$ then $\lambda$ is a Magidor cardinal (by definition, rememeber that $F$ is uniform), and hence satisfies Theorem \ref{t1}. The following lemma shows that one may assume that $\jmath ''\lambda$ belongs to the filter.

\begin{lemma}
\label{lem2} Assume $G$ is a Magidor filter over $\lambda$. \newline 
Then there exists an elementary embedding $\jmath:M\rightarrow {\rm V}_{\lambda+\omega}$ such that ${\rm crit}(\jmath)<\lambda=\jmath(\lambda)$ and $\jmath ''\lambda\in G$.
\end{lemma}

\par\noindent\emph{Proof}. \newline 
Let $G$ be a Magidor filter over $\lambda$. Let $N$ be an elementary submodel of ${\rm V}_{\lambda+\omega}$ of size $\lambda$, so that $\lambda+1\subseteq N$. Fix a well ordering of ${\rm V}_{\lambda+\omega}$ in $N$. Being a Magidor cardinal, $|[\lambda]^{\aleph_0\text{-bd}}|=\lambda$. enumerate all the $\aleph_0$-bounded subsets of $\lambda$ by $\{t_i:i<\lambda\}$ in such a way that if ${\rm sup}(t_i)\leq\alpha$ then $i<|\alpha|^{\aleph_0}$. Denote the map $t_i\mapsto i$ by $t$.

Let $\langle f_n:n\in\omega\rangle$ be a set of Skolem functions for $({\rm V}_{\lambda+\omega},\in,t)$. We clump these functions into one single function $g:[\lambda]^{\aleph_0\text{-bd}}\rightarrow{\rm V}_{\lambda+\omega}$, and then we define $h:[\lambda]^{\aleph_0\text{-bd}}\rightarrow\lambda$ by $h(s)=g(s)$ if $g(s)\in\lambda$ and $h(s)=0$ otherwise.
By the Magidority of $G$ there exists $A\in G$ for which $h''[A]^{\aleph_0\text{-bd}}\neq\lambda$. Notice, however, that $A\subseteq h''[A]^{\aleph_0\text{-bd}}$ since the identity is one of our Skolem functions.

Let $M'=g''[A]^{\aleph_0\text{-bd}}$, so $M'\prec N\prec {\rm V}_{\lambda+\omega}$ and $A\subseteq M'$. Let $M$ be the Mostowski collapse of $M'$ and let $\jmath: M\rightarrow M'$ be the inverse of the collapse. It follows that $\jmath'' \lambda=M'\cap\lambda\supseteq A$ and hence $\jmath'' \lambda\in G$, so we are done.

\hfill \qedref{lem2}

We can state now our main result. The proof is modelled after Tryba, \cite{MR877853}, who showed that there are many J\'onsson cardinals below the first cardinal which carries a J\'onsson filter. Our plan is to prove that there would be many $\omega$-J\'onsson cardinals below the first cardinal which carries a Magidor filter, and this is impossible under the axiom of choice.

However, one key-feature in the proof of Tryba fails when moving to Magidority. The main idea of Tryba is to project a J\'onsson filter with some function $g:\lambda\rightarrow\lambda$, and J\'onssonicity is preserved under taking the preimage of any such function. But the projection of a Magidor filter with a function that keeps its uniformity need not be a Magidor filter. The subtle point is that $\aleph_0$-bounded subsets may transfer into unbounded subsets of size $\aleph_0$. 
In order to cope with this problem, we impose another requirement on the projecting function.
It turns out that we can define a suitable function which keeps Magidority.

\begin{theorem}
\label{mt} The main theorem. \newline 
Assuming the axiom of choice, there are no Magidor filters.
\end{theorem}

\par\noindent\emph{Proof}. \newline 
Assume toward contradiction that $F$ is a Magidor filter over $\lambda$. We may assume that $F$ is an ultrafilter, as any ultrafilter extending $F$ is a Magidor filter as well. By virtue of Lemma \ref{lem1}, choose a regular $\theta<\lambda$ such that $F$ is not $(\omega,\theta)$-regular.

We shall prove that there is a function $g:\lambda\rightarrow\lambda$ such that $|g^{-1}(\{\gamma\})|<\theta$ for every $\gamma<\lambda$ (we say that $g$ is almost one-to-one), $g$ is monotonic and unbounded, and $g$ is $<_F$-minimal with respect to these properties. We call $g$ the projecting function, and we shall use it in order to create a Magidor filter with a kind of weak normality.

Assume, towards a contradiction, that no such $g$ exists. We construct a sequence of functions $\langle f_\varepsilon:\varepsilon<\theta\rangle$, each of which is a function from $\lambda$ into $\lambda$, with the following properties:
\begin{enumerate}
\item [$(a)$] $|f_\varepsilon^{-1}(\{\gamma\})|<\theta$ for any $\gamma<\lambda$ and every $\varepsilon<\theta$.
\item [$(b)$] $f_\varepsilon$ is monotonic and unbounded in $\lambda$ for every $\varepsilon<\theta$.
\item [$(c)$] $\varepsilon<\zeta\Rightarrow f_\zeta<_F f_\varepsilon$.
\item [$(d)$] $\varepsilon<\zeta\Rightarrow f_\zeta\leq f_\varepsilon$.
\end{enumerate}

The construction starts from the identity function on $\lambda$ as $f_0$. 
Clearly, all the requirements are satisfied.
In the successor stage $\varepsilon+1$ we employ the assumption toward contradiction which means that $f_\varepsilon$ does not satisfy the above mentioned properties in order to choose some $h<_F f_\varepsilon$ which is almost one-to-one, monotonic and unbounded. We set $f_{\varepsilon+1} = {\rm min}(h,f_\varepsilon)$. At limit stages we let $f_\varepsilon(\alpha)={\rm min}\{f_\beta(\alpha):\beta<\varepsilon\}$ for every $\alpha<\lambda$. Notice that monotonicity is still preserved.

By $(b)$, let $A_\varepsilon=\{\alpha<\lambda: f_{\varepsilon+1}(\alpha)< f_\varepsilon(\alpha)\}$ for every $\varepsilon<\theta$. The collection $\mathcal{A}=\{A_\varepsilon:\varepsilon<\theta\}$ is a subcollection of $F$. Keep in mind that $F$ is not $(\omega,\theta)$-regular, and choose $\mathcal{B}\subseteq\mathcal{A}, |\mathcal{B}|=\aleph_0$ so that the intersection of the members of $\mathcal{B}$ is not empty. By choosing any ordinal $\gamma$ in this intersection, we have an infinite decreasing sequence of ordinals (think of $f_\varepsilon(\gamma)$ for each member of $\mathcal{B}$), a contradiction.

We use the function $g$ in order to project the ultrafilter $F$ and get a new filter $G$ over $\lambda$ as follows:
$$
G=\{x\subseteq\lambda: g^{-1}[x]\in F\}.
$$
The fact that $G$ is an ultrafilter is immediate. The fact that $|g^{-1}(\{\gamma\})|<\theta$ for every $\gamma<\lambda$ ensures that $G$ is uniform, as if $|x|<\lambda$ then $|g^{-1}[x]|<|x|\cdot\theta<\lambda$ and hence $g^{-1}[x]\notin F$. It follows that $G$ contains all the end-segments of $\lambda$.

We proceed to showing that $G$ is Magidor. Let $\delta<\lambda$ be large enough, and assume $f:[\lambda]^{\aleph_0\text{-bd}}\rightarrow\delta$. We define another coloring $h:[\lambda]^{\aleph_0\text{-bd}}\rightarrow\delta$ by $h(t)=f(g''t)$. Choose an element $y\in F$ so that $h''[y]^{\aleph_0\text{-bd}}\neq\delta$.
Let $x=g''y$, so $y\subseteq g^{-1}[x]$ and hence $g^{-1}[x]\in F$ and $x\in G$. Fix an ordinal $\gamma\in\delta - h''[y]^{\aleph_0\text{-bd}}$. We claim that $\gamma\notin f''[x]^{\aleph_0\text{-bd}}$, thus proving that $G$ is a Magidor filter.

Toward contradiction assume that $\gamma=f(t)$ for some $t=\{t_n:n\in\omega\} \in [x]^{\aleph_0\text{-bd}}$. For every $n\in\omega$ choose $s_n\in y$ such that $g(s_n)=t_n$. Notice that $s=\{s_n:n\in\omega\}$ is bounded in $\lambda$, so $s\in[y]^{\aleph_0\text{-bd}}$. 
Indeed, $t$ is bounded in $\lambda$ so one can choose an ordinal $\tau\in({\rm sup}(t),\lambda)\cap x$. Let $\sigma=g^{-1}(\tau)<\lambda$ and observe that $s\subseteq\sigma$ and hence bounded in $\lambda$.
However, $h(s)=f(g''s)=f(t)=\gamma$, a contradiction.

The last property of $G$ that we need is a very weak version of normality. We shall prove the following general fact: if $h:\lambda\rightarrow\lambda$ is a monotonically increasing and regressive function on a set $x\in G$ then $|h^{-1}(\{\gamma\})|\geq\theta$ for some $\gamma<\lambda$. Indeed, since $h$ and $g$ are monotonically increasing, $h\circ g$ is monotonically increasing as well. Set $y=g^{-1}[x]\in F$, and notice that $\beta\in y\Rightarrow h(g(\beta))<g(\beta)$ as $g(\beta)\in x$. Hence, $h\circ g<_F g$ and by the choice of $g$ there is an ordinal $\gamma<\lambda$ for which $|(h\circ g)^{-1}(\{\gamma\})|\geq\theta$. However, $(h\circ g)^{-1}(\{\gamma\}) = g^{-1}(h^{-1}(\{\gamma\}))$, so $h^{-1}$ must be of size at least $\theta$ for some ordinal since $\theta=\cf(\theta)$.
It follows, in particular, that $g$ is not regressive on a set from $G$.

By Lemma \ref{lem2} we choose an elementary embedding $\jmath:M\rightarrow {\rm V}_{\lambda+\omega}$ such that ${\rm crit}(\jmath)<\lambda=\jmath(\lambda)$ and $\jmath ''\lambda\in G$. Let $A\in G$ be the generating set of the models $M^*_n$ from the proof of Theorem \ref{t1}. 
We may assume that $A\subseteq\jmath'' \lambda$, by intersecting $A$ with $\jmath'' \lambda$ (recall that both sets are elements of $G$).
We define a function $f:A\rightarrow\lambda$ by $f(\alpha)={\rm otp}(A\cap\alpha)$. Observe that $f$ is monotonically increasing, unbounded in $\lambda$, one-to-one and regressive over the set $A'=\{\alpha\in A:{\rm otp}(A\cap\alpha)<\alpha\}$. Consequently, $A'\notin G$ and hence $A\setminus A'\in G$.

Choose an ordinal $\alpha\in A\setminus A'$ so that ${\rm crit}(\jmath)<|\alpha|\leq\alpha$. Choose a large enough $n\in\omega$ for which $\alpha<\mu_n<\lambda$. Let $\pi_n:M^*_n\rightarrow M_n$ be the Mostowski collapse, and let $\jmath_n=\pi^{-1}_n$. Notice that ${\rm crit}(\jmath_n)\leq{\rm crit}(\jmath)<|\alpha|$ and $\jmath_n(\alpha)=\alpha$ as ${\rm otp}(A\cap\alpha)=\alpha$. 
It follows that $\jmath_n(|\alpha|)=|\alpha|$ as well.
A focal point here is that $\jmath_n(|\alpha|^{M_n})=|\alpha|^V$ (by the fact that $\jmath_n$ is increasing).

To sum up, we have an elementary embedding $\jmath_n$ from a transitive set $M_n$ into ${\rm V}_\gamma$ such that ${\rm crit}(\jmath_n)<|\alpha|$ and $\jmath_n(|\alpha|)=|\alpha|$. Since $|\alpha|<\mu_n$ we may assume that $(M_n\cap|\alpha|)^{\aleph_0}\subseteq M_n$, see Remark \ref{rrremark}. But this contradicts Claim \ref{cor1}, so we are done.

\hfill \qedref{mt}

The above proof is illuminating in the following sense. Magidor cardinals are strongly connected with rank-into-rank embeddings. In fact, I1 and I2 are Magidor cardinals. These axioms are located on the verge of inconsistency. The additional feature of having a filter which keeps the Magidority, traverses the border into inconsistency.

Actually, the existence of an $\omega$-J\'onsson cardinal is inconsistent with ZFC, and by the above proof the existence of a Magidor filter is another step. It yields an unbounded set which consists of $\omega$-J\'onsson cardinals. This fact gives an insight into the profound difference between J\'onsson cardinals and Magidor cardinals. For J\'onsson cardinals, the filterhood is a strengthening of the consistency strength but it remains in the realm of consistency (every measurable cardinal carries a J\'onsson filter). For Magidor cardinals, the filterhood sends us into inconsistency.

There is another issue to be mentioned in this context. The results in this paper are proved in the frame of ZFC. Without the axiom of choice, the existence of a Magidor filter is possible. Indeed, under AD the measurable ultrafilter over $\aleph_1$ is a Magidor filter.
We may conclude that Magidor filters can serve as a natural large cardinals axiom above $\omega$-J\'onsson cardinals.
We emphasize, however, that the proof of the main theorem employs the axiom of choice, and we do not know what is the relationship between Magidor cardinals and $\omega$-J\'onsson cardinals without AC. Actually, it is not clear whether there exists a Magidor filter over a singular cardinal under weak versions of choice or even in ZF.

Anyhow, the ZFC proof casts Magidor filters as a germane step above $\omega$-J\'onssonicity. Recall that every normal measure over a measurable cardinal contains a measure-one set of Ramsey cardinals, a normal measure over a supercompact cardinals contains a measure-one set of measurable cardinals, and a Magidor filter (if such existed) would imply the existence of an unbounded set of $\omega$-J\'onsson cardinals and a measure-one set of ordinals which satisfy the defining property of $\omega$-J\'onssonicity.

The following figure demonstrates this idea:

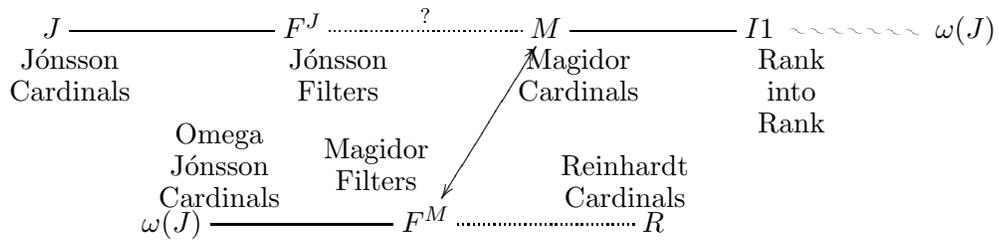
\xymatrix{
J \ar@{-}[rr]_<*+\txt{J\'onsson\\Cardinals} &&F^J \ar@{.}[rr]^-{?}_<<*+\txt{J\'onsson\\Filters} &&M \ar@{-}[rr]_<<*+\txt{Magidor\\Cardinals}  \ar@{<->}[ddl]_>>*+\txt{Magidor\\Filters}
&&I1 \ar@{~}[rr]_<<*+\txt{Rank\\into\\Rank} &&\omega(J) \\
& \\
&\omega(J) \ar@{-}[rr]^<<*+\txt{Omega\\J\'onsson\\Cardinals}
&&F^M \ar@{.}[rr]^>>*+\txt{Reinhardt\\Cardinals} &&R
}

We display a fragment of the chart of large cardinals in two rows.
The second row begins with $\omega$-J\'onssonicity, known to be inconsistent with ZFC but possible under ZF.
The first row reaches up to the point of ZFC-inconsistency, beginning with J\'onsson cardinals. The consistency strength of J\'onsson filters is strictly above the consistency strength of J\'onsson cardinals. For instance, in the canonical model $L[U]$ of Kunen in \cite{MR0277346}, there is only one measurable cardinal, and a measure-one set of J\'onsson cardinals in the unique normal measure. However, only the measurable cardinal carries a J\'onsson filter in $L[U]$. Actually, the consistency strength of a J\'onsson filter over a regular cardinal is measurability (see \cite{MR645907}), so in $L[U]$ we have but one such cardinal.

We do not know whether any Magidor cardinal carries a J\'onsson filter. Of course, natural Magidor cardinals which come from rank-into-rank embedding are limit of measurable cardinals and hence carry a J\'onsson filter. 
We suspect, however, that the consistency strength of Magidority is above measurability (see \cite{GaHa}).
The axiom I1 is strictly stronger than Magidority (again, see \cite{GaHa}), but we do not know if the first Magidor cardinal must be below I2.

The next stage, i.e. $\omega$-J\'onssonicity, is beyond consistency, assuming the axiom of choice. The existence of Magidor filters is another step forward, as shown in this paper. As noted above, Magidor filters may exist without the axiom of choice, e.g. under AD. Reinhardt cardinals point to another open problem, and it is opaque whether their existence can be refuted from the axioms of ZF alone.

We conclude with the following:

\begin{question}
\label{q0} Assume {\rm AD}. \newline 
Does there exist a cardinal $\lambda>\cf(\lambda)=\omega$ which carries a Magidor filter?
\end{question}

\newpage

\bibliographystyle{amsplain}
\bibliography{arlist}

\end{document}